\documentclass[psfig,reqno,a4paper,11pt]{amsart}
\usepackage{amsfonts}

\theoremstyle{plain}
\begingroup

\newtheorem{theorem}{Theorem}[section]
\newtheorem{lemma}[theorem]{Lemma}
\newtheorem{proposition}[theorem]{Proposition}

\newtheorem{remark}[theorem]{Remark}

\newtheorem{example}[theorem]{Example}
\theoremstyle{definition}
\theoremstyle{remark}
\numberwithin{equation}{section}
\newcommand{\var}{\varphi}
\newcommand{\e}{\varepsilon}
\newcommand{\Om}{\Omega}

\newcommand{\dx}{\,dx}

\newcommand{\intcauchy}{\mskip 3mu -\mskip -19mu \int}

\newcommand{\R}{{\mathbb R}}

\newcommand{\salt}{\noalign{\vskip .2truecm}}
\newcommand{\parent}[3]{\left #1 {#3} \right #2} 
\newcommand{\graffe}[1]{\parent \{ \}{#1}} 

\newcommand{\N}{{\mathbb N}}

\newcommand{\bord}{{\partial\Om}}

\input psfig.sty
\newcommand{\bfig}[2]{\begin{figure}\begin{center}\begin{picture}(341.8,#2)(
#1,0)}
\newcommand{\efig}[2]{\end{picture}\caption{#2.}\lbl{#1}\end{center}
\end{figure}}
 
\title[A duality approach for variational problems in domains with cracks]
{A duality approach for variational problems\\ in domains with cracks}
\author[Fran\c{c}ois Ebobisse]{Fran\c{c}ois Ebobisse}
\address[Fran\c{c}ois Ebobisse]{S.I.S.S.A., Via Beirut 2-4, 34014,
Trieste,
Italy and Department of Maths \& applied Maths, University of Cape Town,
 Rondebosch 7700, South Africa}
\email{ebobisse@maths.uct.ac.za}
\author[Marcello Ponsiglione]{Marcello Ponsiglione}
\address[Marcello Ponsiglione]{S.I.S.S.A., Via Beirut 2-4, 34014,
Trieste, Italy}
\email{ponsigli@sissa.it}

\begin{document}
\baselineskip3.3ex

\vskip .2truecm
\begin{abstract}
                
\small{In this paper we study the asymptotic behaviour of the solutions 
of some minimization problems for integral functionals with
 convex integrands, in two-dimensional domains with cracks, 
under perturbations of the cracks in the Hausdorff metric. 
In the first part of the paper, we examine conditions for the stability 
of the minimum problem via duality arguments in convex optimization. 
In the second part, we study the limit problem in some special cases 
when there is no stability, using the tool of $\Gamma$-convergence.
}
\end{abstract}
\maketitle
{\small \tableofcontents}

\section{Introduction}\label{Intro}
Let $\Om $ be a bounded connected and simply connected open subset of $\R ^2$, let $K$ be a compact 
subset of $\overline\Om$ and let $g\in W^{1,p}(\Om)$.  
We consider the following variational problem:
$$
\min_{w=g\,\,on\,\,\partial_D\Om\setminus
K}\int_{\Om\setminus K}f(x,\nabla w)\dx,\leqno(P)
$$
where $\partial_D\Om$ is a non-empty part of the boundary of $\Om$ with a 
finite number of connected components and the function $f:\Om\times\R^2\to\R$ is a Borel function which 
satisfies the following assumptions: there exist 
positive constants $\alpha$, $\beta$, $\gamma$ such 
that, for almost every $x\in\Om$ and for every $\xi\in\R^2$
\begin{eqnarray}
&\alpha|\xi|^p\leq f(x,\xi)\leq\beta |\xi|^p+\gamma;\label{if1}\\
\salt
&f(x,\cdot)\mbox{ is strictly convex}.\label{if2}
\end{eqnarray}

Our purpose in this paper is to study the asymptotic behaviour of the solutions $u_K$ of the problem ($P$)  
with respect to the variations of the compact set $K$ in the Hausdorff metric. 
This problem has been recently studied in \cite{dmro1} for $f(x,\xi)=|\xi|^2/2$ in order 
to give a precise mathematical formulation for the quasi-static growth of brittle fractures, 
following Griffith's criterion of crack growth. 

The study of the asymptotic behaviour of  solutions 
of variational problems with respect to domain variations is also related to some shape optimization 
problems, where very often the nonexistence of solutions is due to the non stability of the state equation.
 By stability  of  problem $(P),$ more  precisely  stability of a given compact set $K$  along
 a sequence $(K_h)$ converging to $K$ in the Hausdorff distance, we mean the convergence in a suitable topology of the sequence of solutions $(u_{K_h})$ of $(P)$ to the
 function $u_K$. It is known that a necessary condition for stability is the convergence of 
the two-dimensional
Lebesgue measure $|K_n|$ of $K_n$ to the two-dimensional
Lebesgue measure $|K|$ of $K$ (see \cite{dmebpo}).

If $f(x,\cdot)$ is differentiable, then the solution $u_K$ solves a nonlinear mixed type boundary value problem. In the literature there are
 various results on the asymptotic behaviour of solutions of elliptic PDE with purely Dirichlet boundary conditions, 
with respect to domains variations. In this case the type of limit problem is known even when there is no stability 
(see for instance \cite{DM}, \cite{budm}, \cite{dmmur}).

Concerning stability results for purely Neumann problems, we can 
mention for instance the papers \cite{chen}, \cite{chdo},  \cite{buzo2},  
\cite{bv1},  \cite{bv2},  \cite{dmebpo}, where the families of 
domains satisfy suitable structural assumptions. In the literature there are
 well known examples showing that without these structural assumptions 
some additional term (typically depending on jumps on the limit set $K$)
 may appear in the limit problem (see \cite{mur}, \cite{daml}, \cite{cor}).  
However, unlike Dirichlet problems, there is not a general characterization 
of the limit problem with Neumann  conditions.

In the first part of this paper we prove the following stability result using 
the duality argument of convex optimization.
\begin{theorem}\label{imainr}
Let $(K_h)$ be a sequence of compact subsets of $\overline\Om$
which converges to a compact set $K$ in the Hausdorff metric. 
Assume that $K_h$ has a uniformly bounded number of connected components,  
 $| K_h|$ converges to $| K|$, and that the intersection of the limits of two different 
connected components of $K_h\cup(\partial\Om\setminus\partial_D\Om)$ is either empty or has
positive $(1,q)$-capacity, where $q$ is the conjugate exponent of $p$. 
Then the compact set $K$ is stable for the problem $(P)$ along the sequence $(K_h)$.
\end{theorem}
When $p\leq 2$ the stability result follows 
immediately from \cite[Theorem 6.3]{dmebpo} even when $\Om$ is not 
simply connected. 
 
The approach by duality consists in proving the stability of the 
limit set $K$ for problem $(P)$
from its stability for the dual problem, which is more easy. Indeed,  
unlike problem $(P)$, 
the admissible functions in the dual problem
 for the approximating sequence $(K_h)$ belong all to the same space 
$W^{1,q}(\Om)$, with the constraint   
 that these functions are constant on every connected components 
of $K_h\cup(\partial\Om\setminus\partial_D\Om)$.
Then the assumptions of Theorem \ref{imainr} give the same constraint 
for the limit set $K$.
 
In the second part of the paper we study several examples of 
non stability in the case $p>2$,
 using the tool of $\Gamma$-convergence.
 For instance,  Example \ref{ex1}  shows that without the capacitary 
assumption in Theorem   
 \ref{imainr}, we may have non stability even when $K_h$ has just two 
connected components.
In the case of non stability, we do not yet have a general 
characterization of the limit problem. However, 
in Example \ref{exfine}, 
we are able to find the limit problem under some geometrical assumptions on the sequence $(K_h)$.

\section{Notation and preliminaries}\label{Notprel} 
Let $\Om $ be a bounded connected and simply connected open subset of $\R^2$   
 with Lipschitz continuous boundary $ \partial\Om$. Let    
 $\partial_D\Om\subset\partial\Om$ be a (non-empty) relatively open subset of  
$\bord$ composed of a finite number of connected components and  
$\partial_N\Om:=\bord\setminus\partial_D\Om$.  
 
Let ${\mathcal K}(\overline\Om)$ 
be the class of compact subsets of $\overline\Om$  
and  ${\mathcal K}_m(\overline\Om)$ be the subset of ${\mathcal K}(\overline\Om)$  
whose elements have at most $m$ connected components.  
 
For any $x\in \Om $ and $\rho>0$, $B_\rho(x)$ denotes  
  the open ball of  $\R^2 $ centered at $x$ with radius $\rho $.  
For any subset $E$ of $\R^2$, $1_E$ is the characteristic function of $E$, $E^c$ is the complement of $E$, and  
$|E|$ is the Lebesgue measure of $E$. Given a subset $F$ of some vectorial space $X$, $I_F$ will denote the indicator function 
 of $F$, i.e., $I_F(x)$ is equal $0$ if $x\in F$ and $+\infty$ otherwise. 
 
Throughout the paper $B$ is an open ball containing $\overline\Om$ and  
$p$ and $q$ are real numbers, with $1<p,\,q<+\infty$ and $p^{-1}+q^{-1}=1$. 
  
\subsection{Conjugate function and duality argument in optimization} 
In this section we recall the concept of duality for the  minimization 
of convex functionals. For more details, the reader is referred to \cite{ET}. 
 
Let $X$ be a reflexive Banach space and let $X^*$ be its topological dual. 
Given a function $F: X\to\overline\R$ convex, lower semicontinuous and proper, 
 the conjugate function 
$F^*: X^*\to\overline\R$ of $F$ is defined by: 
 
$$ 
F^*(u^*):= \sup_{u\in X} \{ \langle u,u^*\rangle - F(u)\}\quad \forall 
u^* \in X^* 
$$  
where $\langle\cdot,\cdot\rangle$ denotes the duality brackets between $X$ and 
 $X^*$.\\ 
We recall that for functionals of the type $F(u)=\int_\Om f(x,u(x)) \dx$, defined in $L^p(\Om,\R^2)$, 
where $f$ satisfies for instance the assumptions (\ref{f1})-(\ref{f2}) below, the following formula holds 
(see for instance \cite[Proposition 2.1]{ET}) 
\begin{equation}\label{conjpart} 
F^*(u^*)=\int_\Om f^*(x,u^*(x)) \dx\qquad \forall u^*\in L^q(\Om,\R^2). 
\end{equation} 
 
Now we consider the following minimization problem 
$$ 
\min_{u\in X} F(u). \leqno (P) 
$$ 
 
Let $Y$ be a Banach space and let $Y^*$ be its topological dual. 
The duality argument in the study of Problem $(P)$ is described as follows. 
We consider a family of perturbations of \mbox{Problem $(P)$:} 
 
$$ 
\min_{u\in X} \Phi(u,\xi) \leqno (P_\xi)  
$$ 
 
where  $\Phi:X\times Y \to \overline\R$ is a convex, lower semicontinuous  
and proper function such that 
 
$$ 
\Phi(u,0) = F(u) \qquad \forall u \in X. 
$$ 
 
The dual problem of $(P)$ with respect to $\Phi$ is given by: 
 
$$ 
\sup_{\xi^*\in Y^*} \{- \Phi^*(0,\xi^*)\}. \leqno (P^*) 
$$ 
The following proposition is proved in \cite[Proposition 2.4]{ET}. 
\begin{proposition}\label{dualita} 
Assume that $\inf_X F$ is finite, that $F$ is coercive and that there exists  
$u_0 \in X$ such that $\xi\to \Phi(u_0,\xi)$ takes values in $\R$ and  
is continuous in $0$. 
 
Then the problems $(P)$ and $(P^*)$ each have at least one solution. Moreover 
\begin{equation}\label{dualform} 
\inf_{u\in X} F(u) = \sup_{\xi^*\in Y^*} \{- \Phi^*(0,\xi^*)\} 
\end{equation} 
 
and the following  relation is satisfied 
 
\begin{equation}\label{reldual} 
\Phi(\bar{u},0) + \Phi^*(0,\bar{\xi}^*) =0, 
\end{equation} 
 
where $\bar{u}$ is a solution of $(P)$ and $\bar{\xi}^*$  
is a solution of $(P^*)$. 
 
Conversely, if $\bar{u}\in X$ and $\bar{\xi}^*\in Y^*$ satisfy 
(\ref{reldual}), then  $\bar{u}$ is a solution of $(P)$ and $\bar{\xi}^*$  
is a solution of $(P^*)$. 
\end{proposition} 
\vskip .2truecm 
In this paper, we will deal with functionals $\Phi$ of this type:  
\begin{equation}\label{fpart} 
\Phi(u,\xi) = F_1(u) + F_2(Au - \xi), 
\end{equation}  
where $F_1 :X\to\R$ 
and $F_2 :Y\to\R$ are convex lower semicontinuous functions and $A:X\to Y$ is a linear continuous operator. 
In this case, we have that 
 
\begin{equation}\label{formulad} 
\Phi^*(0,\xi^*)= F_1^*(A^* \xi^*)+ F^*_2(-\xi^*). 
\end{equation} 
 
where  $A^*:Y^*\to X^*$ denotes the transpose of the operator $A$. 
\subsection{Deny-Lions spaces} 
 Given an open subset 
$U$ of $\R^2$, the Deny-Lions space is defined by 
$$L^{1,p}(U):=\{u\in L^p_{\rm loc}(U):\, \nabla u\in L^p(U,\R^2)\}.$$ 
It is well-known that $L^{1,p}(U)$ coincides with the Sobolev space $W^{1,p}(U)$  
whenever $U$ is bounded and has a Lipschitz continuous boundary.  
It is also known that the set $\{\nabla u:\, u\in L^{1,p}(U)\}$ is a closed subspace of 
$L^p(U,\R^2)$. The Deny-Lions spaces $L^{1,p}$ are usually involved in minimization 
 problems of the type (\ref{eq1}) below  in non-smooth domains, 
 where Poincar\'e inequalities  do not hold in general. 
 For further properties of the spaces $L^{1,p}$ we refer the 
reader to \cite{deli} and \cite{maz1}. 
 
\subsection{The minimization problem}  
Let $f:\Om\times\R^2\to\R$ be a Borel function which  
satisfies the following assumptions: there exist  
positive constants $\alpha$, $\beta$, $\gamma$ such  
that, for almost every $x\in\Om$ and for every $\xi\in\R^2$ 
\begin{eqnarray} 
&\alpha|\xi|^p\leq f(x,\xi)\leq\beta |\xi|^p+\gamma;\label{f1}\\ 
\salt 
&f(x,\cdot)\mbox{ is strictly convex}.\label{f2} 
\end{eqnarray} 
Given $K\in {\mathcal 
K}(\overline\Om)$ and a function $g\in W^{1,p}(\Om)$,  
we consider the following  minimization problem 
\begin{equation}\label{eq1} 
\displaystyle\min_{w}\graffe{\int_{\Om\setminus K}f(x,\nabla w)\dx\,\mbox{: }w\in  
L^{1,p}(\Om\setminus K)\,,\quad w=g\,\mbox{ on }\,\partial_D\Om\setminus K} 
\end{equation} 
whose solution exits from  direct methods of the calculus of variations and is unique in the sense of gradients.

\subsection{$\Gamma$-convergence} Let us recall the definition of  
De Giorgi's {\it $\Gamma$-convergence} in metric spaces. 
 Let $(X,d)$ be a metric space. We say that a sequence $F_h:X\to [-\infty ,+\infty 
  ]$ $\Gamma $-converges to $F:X\to [-\infty ,+\infty ]$ (as $h\to\infty$) if for all $u \in X$ we have 
\begin{itemize} 
\item[{\rm (i)}] ({\it lower limit inequality}) for every sequence $(u_h)$ converging to 
  $u$ in $X$,  
$$ \liminf\limits _{h\to\infty }F_h(u_h)\geq F(u);$$ 
\item[{\rm (ii)}] ({\it existence of a recovery sequence}) there exists a sequence 
  $(u_h)$ converging to $u$ in $X$, such that 
$$\limsup\limits _{h\to \infty }F_h(u_h)\leq F(u).$$ 
\end{itemize} 
The function $F$ is called $\Gamma$-limit of $(F_h)$  
(with respect to $d$), and we write $F=\Gamma\mbox{-}\lim_{h}F_h$.  
The peculiarity of this type of convergence is its variational 
character explained in the following proposition. 
\begin{proposition}\label{Gamma-conv-prop} 
 Assume that $\{F_h\}$ $\Gamma $-converges to $F$ and that there 
  exists a compact set $K\subseteq X$ such that  
$$ \inf\limits _{u\in K}F_h(u)=\inf\limits _{u\in 
  X}F_h(u)\quad\quad\forall h\in \N.$$ 
Then 
\begin{itemize} 
\item[{\rm (i)}] $\inf _XF_h$ converges as $h\to \infty $ to $\min _XF$ and any 
  limit point of any sequence $(u_h)$ such that  
$$\lim\limits _{h\to \infty }\Bigl (F_h(u_h)-\inf\limits _{u\in 
  X}F_h(u)\Bigr )=0$$ 
is a minimizer of $F$. 
\item[{\rm (ii)}] $(F_h+G)$ $\Gamma $-converges to $F+G$ for any 
  $G:X\to ]-\infty ,+\infty [$ continuous. 
\end{itemize} 
\end{proposition} 
We refer the reader to  \cite{dm} for an exhaustive treatment 
of this topic.

\subsection{Hausdorff convergence}\label{haus} 
 The {\it Hausdorff distance}  
between two closed subsets $K_1$ and $K_2$ of $\overline\Om$ is defined by  
$$ 
d_H(K_1,K_2):=\max\graffe{\displaystyle\sup_{x\in K_1}{\rm dist}\,(x,K_2)\,,\, 
\displaystyle\sup_{x\in K_2}{\rm dist}\,(x,K_1)},$$ 
 with the conventions ${\rm dist} 
\,(x,\emptyset)={\rm diam}\,(\Om)$ and $\sup\emptyset =0$, so that 
$$ 
d_H(\emptyset\,,K)= 
\begin{cases} 
0 &\mbox{ if }K=\emptyset, \\ 
{\rm diam}\,(\Om)&\mbox{ if }K\neq\emptyset .  
\end{cases} 
$$  
Let $(K_h)$ be a sequence of compact subsets of $\overline\Om $. We say that $(K_h)$  
converges to $K$ in the {\it Hausdorff metric} if $d_H(K_h\,,K)$ converges to $0$.  
It is well-known (see e.g., \cite[Blaschke's Selection Theorem]{falc}) that  
${\mathcal K}(\overline\Om)$ and ${\mathcal K}_m(\overline\Om)$ are 
 compact with respect to the  Hausdorff convergence. 
 
\vskip .3truecm 
  
In order to study the continuity of the solution $u$ of (\ref{eq1}) with respect to  
the variations of the compact set $K$, we should be able to compare two solutions defined in  
two different domains. This is why, throughout this paper, given a function $u\in L^{1,p}(\Om\setminus K)$,  
 we extend $\nabla u$ in $\Om$ by setting $\nabla u=0$ in $\Om\cap K$.  
 
\subsection{Capacity}\label{cap} 
 Let $1<r<\infty$. We recall that $B$ is a fixed open ball containing  
$\overline\Om$. For every subset $E$ of $B$, the $(1,r)$-capacity of $E$ in $B$, denoted by  
$C_r(E,B)$ or simply by $C_r(E)$ (when there is no ambiguity),  
is defined as the infimum of $ \int_B|\nabla u|^r\,dx$ over the set of all functions  
$u\in W^{1,r}_0(B)$ such that $u\geq 1$ a.e. in a neighborhood of $E$.  
If $r>2$, then $C_r(E)>0$ for every nonempty set $E$. On the  
contrary, if $r=2$ there are nonempty sets $E$ with $C_r(E)=0$ (for  
instance, $C_r(\{x\})=0$ for every $x\in B$). 
 
We say that a property $\mathcal P(x)$ holds $C_r$-{\it quasi everywhere} (abbreviated $C_r$-{\it q.e.}) 
 in a set $E$ if it holds for all $x\in E$ except a subset $N$ of $E$ with $C_r(N)=0$. 
 We recall that the expression {\it almost everywhere} (abbreviated {\it a.e.}) refers, as usual, 
 to the Lebesgue measure.  
 
A function $u:E\to\overline\R$ is said to be {\it quasi-continuous}  
if for every $\e$ there exists $A_\e\subset E$, with  $C_r(A_\e)<\e$, such that the restriction  
of $u$ to $E\setminus A_\e$ is continuous. If $r>2$ every quasi-continuous function is continuous,  
while for $r=2$ there are quasi-continuous functions that are not  
continuous. It is well known that,  
for every open subset $U$ with  $\overline U\subset B$, any function $u\in 
 L^{1,r}(U)$ has a {\it quasi-continuous representative}  
$\overline u:U\cup\partial_LU\to\R$ which satisfies 
$$ 
\lim_{\rho\to 0^+}\intcauchy_{B_\rho (x)\cap U}|u(y)-\overline u(x)|\,dy=0 
\quad\mbox{for $C_r$-q.e. }x\in U\cup\partial_LU, 
$$  
where $\partial_LU$ denotes the Lipschitz part of the boundary $\partial U$ of $U$. 
 We recall that if $u_h$ converges to $u$  strongly in $W^{1,r}(U)$, then  
a subsequence of $\overline u_h$ converges to $\overline u$ pointwise $C_r$-q.e. on $U\cup\partial_LU$. 
 To simplify the notation we shall always identify throughout the paper each function $u\in L^{1,r}(U)$ with 
 its quasi-continuous representative $\overline u$. 
   
For these and other properties on quasi-continuous representatives 
  the reader is referred to  
\cite{EG}, \cite{hekima}, \cite{maz1}, \cite{ziem}. 
\vskip .1truecm 
The following lemma is proved in \cite[Lemma 4.1]{dmro1} for $p=2$.  
The case $p\neq 2$ can be proved in the same way. 
\begin{lemma}\label{lemmaaux} 
Let $(K_h)$ be a sequence in $\mathcal K(\overline\Om)$  
which converges to a compact set $K$ in the Hausdorff  
metric. Let $u_h\in L^{1,p}(\Om\setminus K_h)$ be a  
sequence such that $u_h=0$ $C_p$-q.e. on $\partial_D\Om\setminus K_h$ and 
 $(\nabla u_h)$ is bounded in $L^p(\Om,\R^2)$.  
Then, there exists a function $u\in L^{1,p}(\Om\setminus K)$ with $u=0$ $C_p$-q.e. on  
$\partial_D\Om\setminus K$ such that, up to a subsequence,   
$\nabla u_h$ converges weakly to $\nabla u$ in $L^p(A,\R^2)$  
for every $A\subset\subset\Om\setminus K$.  
If in addition $|K_h|$ converges to $|K|$, then  
 $\nabla u_h$ converges weakly to $\nabla u$ in $L^p(\Om,\R^2)$. 
\end{lemma} 
 
The following lemma will be crucial in the proof of our main results. 
\begin{lemma}\label{moscostante1} 
Let $(K_h)\subset\mathcal K_m(\overline\Om)$ be a sequence  which 
converges to a compact set $K$ in 
the Hausdorff metric, 
and let $(v_h)$ be a sequence in 
$W^{1,q}(\Om)$ which converges to a function $v$ 
weakly in $W^{1,q}(\Om)$. Assume that  the intersection of the limits  
of two different connected components of $K_h$ is either empty or has 
positive $C_q$-capacity and that every function $v_h$ is constant $C_q$-q.e.\ in each connected 
component of $K_h$. 
Then $v$ is constant $C_q$-q.e.\ in each connected 
component of $K$. 
\end{lemma} 
\begin{proof} By extending both functions $v_h$ and $v$ in the open ball $B$ 
containing $\overline\Om$ such that 
$v_h\rightharpoonup v$ weakly in $W^{1,q}(B)$ and arguing as in \cite[Lemma 3.5]{dmebpo} we obtain that $v$ is   
constant $C_q$-q.e.\ in the limit of each connected component of $K_h$.  
Now using the assumption that the intersection of the limits  
of two different connected components of $K_h$ is either empty or has 
positive $C_q$-capacity, we get that $v$ is constant $C_q$-q.e.\ in each connected 
component of $K$. 
\end{proof} 
 
The following lemma will be used in order  
to get the strong convergence of solutions in our main results. 
 
\begin{lemma}\label{stroco} 
Let $f:\Om\times \R^2 \to R$ be a Borel function which satisfies the assumptions 
(\ref{f1}) and (\ref{f2}), and let $(\xi_h)$ be a sequence in $L^p(\Om,\R^2)$  weakly 
converging to some $\xi \in L^p(\Om,\R^2)$.  
If $\int_\Om f\big(x,\xi_h(x)\big) \dx$ converges to  
$\int_\Om f\big(x,\xi(x)\big) \dx$, 
 then $\big(\xi_h\big)$ converges  to $\xi$ strongly in $L^p(\Om,\R^2)$. 
\end{lemma} 
\begin{proof} 
By the convexity of $f$, we have the following lower semicontinuity inequality
$$\int_\Om f\big(x,\xi(x)\big) \dx \leq  \liminf_h   \int_\Om
f\big(x,(\xi+\xi_h(x))/2\big) \dx.$$ Hence,   \begin{eqnarray}\label{sc} 
\nonumber \limsup_h \int_\Om\left[\frac{1}{2} f\big(x,\xi(x)\big)+  
\frac{1}{2}f\big(x,\xi_h(x)\big) - f\big(x,(\xi(x)+\xi_h(x))/2\big)  
\right]\dx \leq \\ 
\leq 
\limsup_h\left[  
\frac{1}{2} \int_\Om f\big(x,\xi_h(x)\big) \dx -  
\frac{1}{2} \int_\Om f\big(x,\xi(x)\big) \dx\right] =0 
\end{eqnarray} 
On the other hand, by the convexity of $f(x,\cdot)$ we have that  
$$ \frac{1}{2}  f\big(x,\xi(x)\big)  + \frac{1}{2}  f\big(x,\xi_h(x)\big) 
-  f\big(x,(\xi(x)+\xi_h(x))/2\big)$$ 
is non negative, and thus 
\begin{equation}\label{sto0} 
 \frac{1}{2}  f\big(x,\xi(x)\big)  + \frac{1}{2}  f\big(x,\xi_h(x)\big) 
-  f\big(x,(\xi(x)+\xi_h(x))/2\big) \to 0 \quad \mbox{ strongly in } L^1(\Om). 
\end{equation} 
Up to a subsequence, we have 
$$ 
 \frac{1}{2}  f\big(x,\xi(x)\big)  + \frac{1}{2}  f\big(x,\xi_h(x)\big) 
-  f\big(x,(\xi(x)+\xi_h(x))/2\big) \to 0 \quad \mbox{ a.e. in } \Om. 
$$ 
By the strict convexity of $f(x,\cdot)$, it easily follows that  
$$\xi_h(x) \to \xi (x) \quad \mbox{ a.e. in } \Om$$ 
and hence $f(x,\xi_h(x))\to f(x,\xi(x))\mbox{ a.e. in }\Om.$ 
Then by Fatou's Lemma we get 
\begin{eqnarray}
\hskip -.5truecm\liminf_h \int_\Om\left[f\big(x,\xi_h(x)\big) + f\big(x, \xi(x)\big)
-|f\big(x,\xi_h(x)\big)-f\big(x,\xi(x)\big)|\right]dx  \geq
2 \int_\Om f\big(x, \xi(x)\big)dx 
\end{eqnarray} 
from which it follows that $\limsup_h \int_\Om 
|f\big(x,\xi_h(x)\big)-f\big(x,\xi(x)\big)|\dx \leq 0$, that is 
\begin{equation}\label{fxi} 
f\big(x,\xi_h(x)\big)\to f\big(x,\xi(x)\big)\mbox{ strongly in
}L^1(\Om).\end{equation}  Now from (\ref{fxi}) and by assumption 
(\ref{f1}), we have (up to a subsequence) that  $(\xi_h)$  
is dominated 
in $L^p(\Om,\R^2)$, which  together with the pointwise convergence above 
 imply that $\xi_h \to \xi$ strongly in $L^p(\Om,\R^2)$. 
\end{proof} 
\section{The dual problem}\label{dualpb} 
 
According to the notation of Section $2.1$, we set 
\begin{eqnarray*}  
&&X:= L^{1,p}(\Om\setminus K),\quad Y:= L^{p}(\Om\setminus K,\R^2);\\ 
\salt 
&&F_1(u): = I_{\{w\mbox{: }w=g\mbox{ on }\partial_D \Om\setminus K\}}(u) \quad \forall u\in L^{1,p}(\Om\setminus K);\\ 
\salt 
&&F_2(\xi): = \int_{\Om\setminus K} f(x,\xi) \dx \quad \forall \xi \in L^{p}(\Om\setminus K,\R^2);\\ 
\salt 
&&A u: = \nabla u\quad \forall u\in L^{1,p}(\Om\setminus K). 
\end{eqnarray*}  
So, the functional to minimize in (\ref{eq1}) is of the type (\ref{fpart}), that is 
$$F(u) = F_1(u) + F_2(A u).$$ 
According to  formula (\ref{formulad}), we need to compute $F_1^*$ and $F_2^*$. 
First of all, for every $u^* \in (L^{1,p}(\Om\setminus K))^*$  
there exists some $\xi^*\in L^q(\Om\setminus K,\R^2)$ such that 
\begin{equation}\label{id} 
\langle u^*,u\rangle = \int_{\Om\setminus K} \xi^* \nabla u \dx \quad \forall u\in L^{1,p}(\Om\setminus K). 
\end{equation} 
Note that 
\begin{equation}\label{a*} 
A^* \xi^*=u^* 
\end{equation} 
Using this representation, we have that 
$$ 
F_1^*(u^*) = \sup_{\substack{u \in L^{1,p}(\Om\setminus K)\\ 
\salt 
u=g\mbox{ on }\partial_D \Om\setminus K}}\,\, 
\Bigl[\,\int_{\Om\setminus K} \xi^* \nabla u \dx \,\Bigr] 
 = \sup_{\substack{u \in L^{1,p}(\Om\setminus K)\\ 
\salt 
u=0\mbox{ on }\partial_D \Om\setminus K}}\,\,\Bigl[\,\int_{\Om\setminus K} \xi^* \nabla u \dx 
+\int_{\Om\setminus K} \xi^* \nabla g \dx\,\Bigr]. 
$$ 
So, by the fact that the supremum of an affine function on a vector space is equal to $0$ or to 
$\infty$, we obtain 
\begin{equation}\label{formulad1} 
F_1^*(u^*) = F_1^*(A^*\xi^*) =  
\begin{cases} 
\displaystyle 
\int_{\Om\setminus K} \xi^* \nabla g \dx & \mbox{ if }  
\begin{cases} 
\int_{\Om\setminus K} \xi^* \nabla \varphi \dx =0 & 
\forall \varphi\in  L^{1,p}(\Om\setminus K),\\  
\varphi=0  \mbox{ on } \partial_D\Om\setminus K, & { } 
\end{cases} 
\\ 
\salt 
+ \infty & \mbox{ otherwise }. 
\end{cases} 
\end{equation} 
Note that the condition $\int_{\Om\setminus K} \xi^* \nabla \varphi \dx =0$  
 $\forall\varphi\in  L^{1,p}(\Om\setminus K)$ with $\varphi=0   
\mbox{ on } \partial_D\Om\setminus K$, is the weak formulation of 
$$\begin{cases} 
{\rm div}\,\xi^*=0 & \mbox{ in }\Om\setminus K,\\ 
\salt 
\xi^*\cdot\nu=0 & \mbox{ on }\partial_N\Om\cup\partial K. 
\end{cases} 
$$ 
On the other hand, from (\ref{conjpart}) we have also 
\begin{equation}\label{formulad2} 
F_2^*(\xi^*)=\int_{\Om\setminus K} f^*(x,\xi^*)\dx\quad \forall \xi^*\in L^q(\Om\setminus K,\R^2). 
\end{equation} 
 
Finally, formula (\ref{formulad}) in this case gives 
 
\begin{equation}\label{formulad3} 
\Phi^*(0,\xi^*) = 
\begin{cases} 
\displaystyle 
\int_{\Om\setminus K}\bigl[ f^*(x,-\xi^*) +\xi^* \nabla g\bigr] \dx & \mbox{ if }  
\begin{cases} 
\int_{\Om\setminus K} \xi^* \nabla \varphi \dx =0 & 
\forall \varphi\in  L^{1,p}(\Om\setminus K),\\ \varphi=0  \mbox{ on } \partial_D\Om\setminus K,&{ } 
\end{cases} 
\\ 
\salt 
+ \infty & \mbox{ otherwise}. 
\end{cases} 
\end{equation} 
 
The duality formula (\ref{dualform}) in this case is given by 
\begin{equation}\label{dualproblem2} 
\min_{u\in L^{1,p}(\Om\setminus K)} F(u) =  
\sup_{\substack{\xi^* \in L^{q}(\Om\setminus K,\R^2)\\ 
\salt 
\substack{\int_{\Om\setminus K} \xi^* \nabla \varphi \dx =0\\ \salt 
\forall \varphi \in L^{1,p}(\Om\setminus K),\,\, 
\varphi=0 \mbox{ on }\partial_D\Om\setminus K}}}  
\hskip-1truecm \int_{\Om\setminus K} \bigl[-f^*(x,-\xi^*) 
-\xi^* \nabla g\bigr] \dx  
\, . 
\end{equation} 
Note that all the results above are actually valid in any dimension,  
while in two dimensional domains, 
 the dual problem in the right hand-side of (\ref{dualproblem2}) 
can be rewritten as a maximum problem in some suitable subspace  
of $W^{1,q}(\Om)$.  
 
To this aim, let $R$ be the rotation on $\R^2$  
defined by $$R(y_1,y_2):=(-y_2,y_1)$$ and let  
$i:W^{1,q}(\Om) \to L^q(\Om,\R^2)$ be the mapping defined by 
\begin{equation}\label{i} 
i(v) := R \nabla v\quad \forall v\in W^{1,q}(\Om). 
\end{equation} 
For every compact set $K\subset \overline\Om$ we set 
$$  W_K^{1,q}(\Om):=\left\{ 
v\in W^{1,q}(\Om),\, \int_\Om v\dx=0 \mbox{ and } 
 v \mbox{ is constant $C_q-$q.e. in 
every $\mathcal C.\mathcal C.$ of } K 
\right\}$$ 
where the notation $\mathcal C.\mathcal C.$ means connected component.\\ 
The following proposition establishes a bijection between the subspace $W^{1,q}_{K\cup \partial_N\Om}(\Om)$ 
and the set of admissible functions for the dual problem in the right hand-side of (\ref{dualproblem2}). 
 
\begin{proposition}\label{bibesu} 
Assume that the compact set $K$ has a finite number of connected components. Then the mapping $i$ defined in (\ref{i}) 
establishes a bijection between the subspace  
$W^{1,q}_{K\cup \partial_N\Om}(\Om)$ and the set of functions 
$$\left\{\xi^* \in L^{q}(\Om\setminus K,\R^2): \,\,  
\int_{\Om\setminus K} \xi^* \nabla \varphi \dx =0 
\quad \forall \varphi \in L^{1,p}(\Om\setminus K),\,\,  
\varphi=0 \mbox{ on }\partial_D\Om\setminus K\right\}.$$ 
\end{proposition} 
\begin{proof} 
Let $v \in W^{1,q}_{K\cup \partial_N\Om}(\Om)$.  
Let $C^1,\ldots,C^l$ be the connected components of  
 $K\cup\partial_N\Om$. Since $v=c^i$ $C_q$-q.e on $C^i$, by  
\cite[Theorem 4.5]{hekima} we can approximate $v$ strongly in $W^{1,q}(\Om)$  
by a sequence of functions $v_n\in C_c^\infty(\R^2)$ that are constant in a  
suitable neighborhood $V_n^i$ of $C^i$. Let $\var\in L^{1,p}(\Om\setminus K)$ with 
 $\var=0$ on $\partial_D\Om\setminus K$ and let $\var_n\in W^{1,p}_0(\Om\setminus K)$  
such that $\var_n=\var$ in $\Om\setminus\bigcup_iV^i_n$. Then we have that  
\begin{equation}\label{eq31} 
\int_{\Om\setminus K}R\nabla v_n\nabla\var\dx=\int_{\Om\setminus K}R\nabla v_n\nabla\var_ndx=0, 
\end{equation} 
where the last equality follows from the fact that the vector field  
$R\nabla v_n$ is divergence free. Then passing to the limit in (\ref{eq31}) for $n\to\infty$, we get  
$$\int_{\Om\setminus K}R\nabla v\nabla\var\dx=0 
\quad\forall\var\in L^{1,p}(\Om\setminus K)\mbox{ with }\var=0\mbox{ on }\partial_D\Om\setminus K.$$ 
So, $i$ maps the space $W^{1,q}_{K\cup \partial_N\Om}(\Om)$ in the set of admissible function in the 
dual problem.\\ 
Now, let $\xi^* \in L^q(\Om\setminus K,\R^2)$ be such that  
$ \int_{\Om\setminus K} \xi^* \nabla \varphi \dx =0\quad \forall \varphi \in L^{1,p}(\Om\setminus K),\,\,  
\varphi=0$ on $\partial_D\Om\setminus K.$ 
By extending $\xi^*$ by zero on $K$ and still denoting this extension by $\xi^*$, we obtain  
$$ \int_{\Om} \xi^* \nabla \varphi \dx =0\quad \forall \varphi \in \mathcal{D}(\Om),$$ 
i.e., ${\rm div}\,\xi^* = {\rm curl}\,(R\xi^*)= 0$ in $ \mathcal{D}'(\Om)$. 
Since $\Om$ is simply connected, there exists $v\in W^{1,q}(\Om)$ such that $R\xi^*=\nabla v$ 
 a.e. in $\Om$. It is not restrictive to assume that $\int_\Om v \dx = 0$.  
So, we have to prove that $v$ is constant on every connected component of $K\cup \partial_N\Om$. 
Let $C$ a connected component of  $K\cup \partial_N\Om$ and, for every $\e>0$, let  
$$C_\e:= \{x\in \overline{\Om} \mbox{: } {\rm dist}(x,C) < \e\}\quad K_\e:= 
K\cup\partial_N\Om\cup \overline{C_\e}.$$ 
Let $\xi^*_\e$ be the solution of the dual problem in the right 
hand of (\ref{dualproblem2}) with $K$ replaced by $K_\epsilon$. 
By the monotonicity of $\Om\setminus K_\e$, it is easy to see that $\xi^*_\e \to \xi^*$  
strongly in $L^q(\Om,\R^2)$ when $\e \to 0$. As above, let $v_\e\in W^{1,q}(\Om)$  
be such that $\int_\Om v_\epsilon \dx = 0$ and $R\xi_\e^*=\nabla v_\e$ a.e. in $\Om$. 
By the fact that $\nabla  v_\e = 0$ in $C_\e,\,\,v_\e \to v$ strongly in $W^{1,q}(\Om)$, and 
that $C\subset\subset C_\e$, 
we get that $v$ is constant  $C_q$-q.e. on $C$, so $v\in W^{1,q}_{K\cup \partial_N\Om}(\Om)$. 
\end{proof} 
 
Using Proposition \ref{bibesu}, the dual problem can be rewritten as 
 
\begin{equation}\label{eqdual} 
\sup_{v\in  W^{1,q}_{K\cup \partial_N\Om}(\Om)} \int_{\Om\setminus K} 
\bigl[- f^*(x,R\nabla v)  
+ R \nabla v \nabla g \bigr]\dx\, . 
\end{equation} 
 
So, the duality formula  (\ref{dualproblem2}) in two dimensional domains 
 becomes  
\begin{equation}\label{dualproblem3} 
\min_{u\in L^{1,p}(\Om\setminus K)} F(u) =  
\sup_{v\in  W^{1,q}_{K\cup \partial_N\Om}(\Om)} \int_{\Om\setminus K} 
\bigl[- f^*(x,R\nabla v) + R \nabla v \nabla g \bigr]\dx\, . 
\end{equation} 
 
Let $u\in L^{1,p}(\Om\setminus K)$ be a solution of the left hand-side of (\ref{dualproblem3}). 
 A solution $v\in W^{1,q}_{K\cup \partial_N\Om}(\Om)$ of the right hand-side of (\ref{dualproblem3}) 
is called {\it conjugate} of $u$. 
 
The duality relation between $u$ and $v$ is 
\begin{equation}\label{dualrel3} 
\int_{\Om\setminus K} f(x,\nabla u) \dx =  
\int_{\Om\setminus K}\bigl[- f^*(x,R\nabla v)+ R \nabla v \nabla g \bigr]\dx. 
\end{equation} 
 
\begin{remark}\label{C1} 
{\rm Since $u=g$ on $\partial_D\Om\setminus K$ and $v\in W^{1,q}_{K\cup \partial_N\Om}(\Om)$ it follows  
that  
$$\int_{\Om\setminus K} R \nabla v \nabla g \dx=\int_{\Om\setminus K} R \nabla v \nabla u \dx.$$ 
Hence (\ref{dualrel3}) becomes 
\begin{equation}\label{dualrel4} 
\int_{\Om\setminus K}\bigl[f(x,\nabla u)+f^*(x,R\nabla v)-R \nabla v \nabla u\bigr]\dx=0. 
\end{equation} 
Since the integrand in (\ref{dualrel4}) is positive, we get 
$$f(x,\nabla u)+f^*(x,R\nabla v)-R \nabla v \nabla u=0\quad\mbox{a.e. in }\Om\setminus K.$$ 
That is  
$$R \nabla v\in\partial_\xi f(x,\nabla u)$$ 
where $\partial_\xi f(x,\nabla u)$ denotes the subdifferential of $f(x,\cdot)$ at the point $\nabla u$.  
Whenever $f(x,\cdot)$ is also of class $C^1$, then $f^*(x,\cdot)$ is strictly convex and hence the dual problem as a unique  
solution $v$ such that  $R \nabla v=\nabla_\xi f(x,\nabla u)$. 
 
For $f(\xi):=\frac{1}{p}|\xi|^p$ we obtain 
$$R \nabla v=|\nabla u|^{p-2}\nabla u\quad\mbox{a.e. in }\Om\setminus K.$$ 
In particular for $p=2$ we obtain the classical notion of harmonic conjugate. 
}\end{remark} 
\section{Stability for the minimum problem} 
Let $(K_h)\subset\mathcal K_m(\overline\Om)$ be a sequence  which 
converges to a compact set $K$ in 
the Hausdorff metric. 
We say that  $K$ is stable for the problem $(P)$ along 
the sequence $(K_h)$ if for every function $f$ that satisfies conditions (\ref{f1})-(\ref{f2}) 
and for every $g\in W^{1,p}(\Om),$ we have  
$$\nabla u_h \to \nabla u \quad \mbox{ strongly in }L^{p}(\Om,\R^2),$$ 
where $u_h$ and $u$ are 
solutions of (\ref{eq1}) in $\Om\setminus K_h$ and in $\Om\setminus K$  
respectively. 
 
In the following theorem, we prove the equivalence  
between the stability of $K$ for the minimum problem (\ref{eq1}) and  for 
 its dual under the condition that $f(x,\cdot)$ is of class $C^1$. 
 
\begin{theorem}\label{equistab} 
Let $(K_h)\subset {\mathcal K}_m(\overline\Om)$ be a sequence 
 which converges to a compact set $K$ in the Hausdorff metric and such that  
 $| K_h|$ converges to $| K|$. Assume that $f(x,\cdot)$ is of class $C^1$.  
Let $g\in W^{1,p}(\Om)$.  
Let $u_h$ and $u$ be 
 solutions of (\ref{eq1}) in  $\Om\setminus K_h$ and in $\Om\setminus K$ 
 respectively. 
Let $v_h$ and $v$ be the solutions of the problem (\ref{eqdual}) 
in  $\Om\setminus K_h$ and  in $\Om\setminus K$ respectively.   
Then $$ \nabla u_h \to \nabla u\mbox{ strongly in }L^p(\Om,\R^2)\quad\mbox{if and only 
if}\quad\nabla v_h\to\nabla v\mbox{ strongly in }L^q(\Om,\R^2).$$ 
\end{theorem} 
\begin{proof} 
Assume that $\nabla v_h \to \nabla v$ strongly in $L^q(\Om,\R^2)$. By (\ref{dualrel3}), we have 
 
\begin{equation}\label{dualrel*} 
\int_{\Om\setminus K_h} f(x,\nabla u_h) \dx = \int_{\Om\setminus K_h} 
\bigl[- f^*(x,R\nabla v_h)+ R \nabla v_h \nabla g \bigr]\dx. 
\end{equation} 
By the growth assumptions (\ref{f1}) on the function $f$, we have that $\nabla u_h$ is bounded in $L^p(\Om,\R^2)$.  
So applying Lemma \ref{lemmaaux} to $u_h-g$, we obtain that $\nabla u_h$ converges (up to a  
subsequence) to $\nabla \tilde u$ weakly in $L^p(\Om,\R^2)$ for some function $\tilde u\in L^{1,p}(\Om\setminus K)$  
 with $\tilde u =g$ on $\partial_D\Om\setminus K$. So passing to the limit in (\ref{dualrel*}) we get 
\begin{eqnarray}\label{dualrel5} 
\nonumber & & \int_{\Om\setminus K} f(x,\nabla\tilde u) \dx  
\leq \liminf_{h\to\infty}\int_{\Om\setminus K_h} f(x,\nabla u_h) \dx 
\leq \limsup_{h\to\infty}\int_{\Om\setminus K_h} f(x,\nabla u_h) \dx \\ 
\nonumber & & = 
\lim_{h\to\infty}\int_{\Om\setminus K_h} \bigl[-f^*(x,R\nabla v_h) +R \nabla
v_h \nabla g\bigr] \dx =   \int_{\Om\setminus K}\bigl[- f^*(x,R\nabla v)+ R
\nabla v\nabla g \bigr]\dx\\   & & = \int_{\Om\setminus K} f(x,\nabla u) \dx 
\end{eqnarray} 
where the last equality follows from the duality relation between $u$ and $v$.  
From (\ref{dualrel5}) and the fact that $f(x,\cdot)$ is strictly 
 convex, we get $\nabla u=\nabla \tilde u$ a.e. in $\Om$, and then 
all the inequalities in (\ref{dualrel5}) are equalities. Therefore,  
all the sequence $(\nabla u_h)$ converges weakly in $L^p(\Om,\R^2)$ and 
$$ 
\lim_{h\to\infty}\int_{\Om\setminus K_h} f(x,\nabla u_h) \dx  
= \int_{\Om\setminus K} f(x,\nabla u)\dx. 
$$ 
Using the convention $\nabla u_h = 0$ on $K_h$, $\nabla u = 0$ on $K$, 
and the fact that $|K_h| \to |K|$, we get also 
\begin{equation}\label{eneconv} 
\lim_{h\to\infty}\int_{\Om} f(x,\nabla u_h) \dx  
= \int_{\Om} f(x,\nabla u)\dx. 
\end{equation} 
Now using the strict convexity of $f(x,\cdot)$, we get from (\ref{eneconv}) and from 
Lemma \ref{stroco} that  $(\nabla u_h)$  converges to $\nabla u$  
strongly in  $L^p(\Om,\R^2)$.\\ 
Viceversa, suppose that  $\nabla u_h \to \nabla u$ strongly in $L^p(\Om,\R^2)$.  
Since $f(x,\cdot)$ is of class $C^1$, from Remark \ref{C1}, we have 
 $R\nabla v_h=f_\xi(x,\nabla u_h)$ a.e. in $\Om\setminus K_h$ and  
$R\nabla v=f_\xi(x,\nabla u)$ a.e. in $\Om\setminus K$. Then from the growth assumptions  
 on $f$ we obtain that $\nabla v_h \to \nabla v$ strongly in $L^q(\Om,\R^2)$. 
\end{proof} 
In the following theorem, we give sufficient conditions  
on the sequence $(K_h)$ which guarantee the stability for 
 Problem (\ref{eq1}). 
\begin{theorem}\label{mainr} 
Let $(K_h)\subset {\mathcal K}_m(\overline\Om)$ be a sequence 
 which converges to a compact set $K$ in the Hausdorff metric and such that  
 $| K_h|$ converges to $| K|$. 
Let $g\in W^{1,p}(\Om)$. 
Let $u_h$ and $u$ be 
 solutions of (\ref{eq1})  in $\Om\setminus K_h$ and in $\Om\setminus K$ respectively. 
Assume that  the intersection of the limits of two different  
connected components of $K_h\cup\partial_N\Om$ is either empty or has 
positive $(1,q)$-capacity. 
Then  $\nabla u_h$ converges strongly to $\nabla u$ in $L^p(\Om,\R^2)$.  
\end{theorem} 
\begin{proof} 
Let $v_h \in W_{K_h\cup\partial_N\Om} ^{1,q}(\Om) $ and  
$v \in W_{K\cup\partial_N\Om} ^{1,q}(\Om)$ be  conjugates of $u_h$ and $u$ respectively.  
Up to a subsequence, $\nabla v_h \rightharpoonup \nabla \tilde{v}$ weakly in 
$W^{1,q}(\Om)$ for some $\tilde{v} \in W^{1,q}(\Om)$.  
By the fact that the intersection of the limits of two different connected components of $K_h\cup\partial_N\Om$ is either empty or has 
positive $(1,q)$-capacity, it follows from Lemma \ref{moscostante1} that $\tilde{v} \in W_{K\cup \partial_N\Om}^{1,q}(\Om)$. 
 
By the growth assumptions (\ref{f1}) on the function $f$, we have that $\nabla u_h$ is bounded in $L^p(\Om,\R^2)$.  
So applying Lemma \ref{lemmaaux} to $u_h-g$, we obtain that $\nabla u_h$ converges (up to a  
subsequence) to $\nabla \tilde u$ weakly in $L^p(\Om,\R^2)$ for some function $\tilde u\in L^{1,p}(\Om\setminus K)$  
 with $\tilde u =g$ on $\partial_D\Om\setminus K$.\\ 
Since $ v\in W_{K\cup \partial_N \Om}^{1,q}(\Om)$, using  \cite[Theorem 4.5]{hekima} we can approximate strongly in $W^{1,q}(\Om)$ the function $ v$ 
with smooth functions $w_n$ which are constant in a suitable neighborhood of any connected component of $K\cup \partial_N\Om,$ and hence constant in 
any connected component of   $K_h\cup \partial_N\Om$ for $h$ big enough. So there exists a subsequence of integers $(h_n)$ such that 
$w_n\in W_{K_{h_n}\cup \partial_N \Om}^{1,q}(\Om)$ and $w_n$ converges strongly in $W^{1,q}(\Om)$ to the function $ v$ as $n\to\infty$. 
 Therefore,  
\begin{eqnarray*} 
\nonumber & & \int_{\Om\setminus K}\bigl[- f^*(x,R\nabla v)+  
R \nabla v\nabla g\bigr] \dx  
= \lim_{n\to\infty}\int_{\Om\setminus K_{h_n}}\bigl[- f^*(x,R\nabla w_n)+  
R \nabla w_n\nabla g\bigr] \dx\\ 
\salt 
\nonumber & &\leq \limsup_{n\to\infty}\int_{\Om\setminus K_{h_n}} 
\bigl[- f^*(x,R\nabla v_{h_n})+ R \nabla v_{h_n}\nabla g\bigr] \dx 
\leq  \int_{\Om\setminus K}\bigl[- f^*(x,R\nabla \tilde v)+  
R \nabla \tilde v\nabla g\bigr]\dx . 
\end{eqnarray*} 
Therefore, since $v$ is a maximizer of the dual problem in $\Om\setminus K$,  
all inequalities in the previous formula are equalities, so 
we obtain 
$$ 
\limsup_{n\to\infty}\int_{\Om\setminus K_{h_n}} 
\bigl[- f^*(x,R\nabla v_{h_n})+ R \nabla v_{h_n}\nabla g\bigr] \dx 
=  \int_{\Om\setminus K}\bigl[- f^*(x,R\nabla  v)+  
R \nabla v\nabla g\bigr] \dx. 
$$ 
Now, using the duality relations between the functions $v$, $u$ 
 on one hand, and $v_{h_n}$, $u_{h_n}$ on the other hand, and then 
 passing to the limit, we obtain 
\begin{eqnarray*} 
\nonumber & & \int_{\Om\setminus K}f(x,\nabla u)\dx\,\, = \,\, \int_{\Om\setminus K} 
\bigl[- f^*(x,R\nabla  v)+ R \nabla  v\nabla g\bigr] \dx\\ 
\salt 
\nonumber &= & 
\limsup_{n\to\infty}\int_{\Om\setminus K_{h_n}} 
\bigl[- f^*(x,R\nabla v_{h_n})+ R \nabla v_{h_n}\nabla g\bigr] \dx 
\,\, =\,\,  \limsup_{n\to\infty}\int_{\Om\setminus K_{h_n}}f(x,\nabla u_{h_n})\dx\\ 
\salt 
& \geq &  
\liminf_{n\to\infty}\int_{\Om\setminus K_{h_n}}f(x,\nabla u_{h_n})\dx\,\,\geq\,\, 
 \int_{\Om\setminus K}f(x,\nabla\tilde u)\dx. 
\end{eqnarray*} 
Since $f(x,\cdot)$ is strictly convex, we get that $\nabla\tilde u=\nabla u$ a.e. in $\Om$. 
Therefore, all the sequence $(\nabla u_h)$ converges weakly in  $L^p(\Om,\R^2)$ to $\nabla u$. 
Now, using that $\int_\Om f(x,\nabla u_h) \dx \to \int_\Om f(x,\nabla u) \dx$, by Lemma \ref{stroco} 
we get that  $(\nabla u_h)$ converges to $\nabla u$ strongly in  $L^p(\Om,\R^2)$. 
\end{proof} 
\begin{remark}\label{caseC1} 
{\rm   
When $f(x,\cdot)$ is of class $C^1$, 
Theorem \ref{mainr} is a consequence of Theorem \ref{equistab}. 
Indeed, the assumption that 
the intersection of the limits of two different  
connected components of $K_h\cup\partial_N\Om$ is either empty or has 
positive $(1,q)$-capacity easily guarantees  
the stability for the dual problem, and hence using Theorem \ref{equistab}  
the stability for Problem (\ref{eq1}) also follows.  
} 
\end{remark} 

\section{Some examples of non stability}

In this section we study some examples for $f(x,\xi)=\frac{1}{p}|\xi|^p$. 
Throughout the section
 we assume that $p>2$.
\subsection{Limit problem via $\Gamma$-convergence}

In the following example, the assumptions of Theorem \ref{mainr} hold. 
We show in this case that the stability result follows also by  
$\Gamma$-convergence arguments.

\begin{center}
\psfig{figure=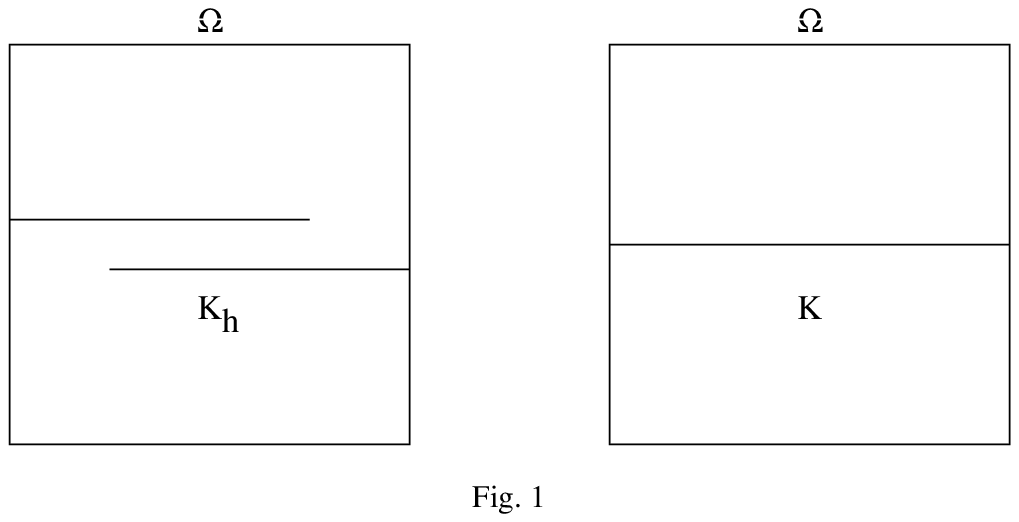}
\end{center}

\begin{example}\label{ex2}
{\rm Let $\Om \, :=\, (-1,1) \times (-1,1)$, $\partial_D \Om := (-1,1)\times \{ -1 , 1\},$\, $K\,=\, [-1,1] \times \{0\}$
 and let 
$$K_h :=\left[-1, \frac{1}{2}\right]\times \left\{\frac{1}{h}\right\} \cup \left[-\frac{1}{2},1\right]\times \left\{-\frac{1}{h}\right\},$$
(see Fig. 1). We consider the sequence of functionals $F_h$ defined in $L^p(\Om)$ by: 
\begin{equation}\label{gammasuc}
F_h(u):= \begin{cases}
\displaystyle\frac{1}{p} \int_{\Om\setminus K_h} |\nabla u|^pdxdy \ & \mbox{ if } u \in W^{1,p}(\Om \setminus K_h)\mbox{ and } u=g \mbox{ on }\partial_D \Om,\\ 
\salt
\displaystyle +\infty & \mbox{ otherwise}.
 \end{cases}\end{equation}

Then, $F_h \, \Gamma$-converges to $F_\infty$ in the strong topology of $L^p(\Om)$, where 
$$F_\infty (u):= \begin{cases}
\displaystyle\frac{1}{p} \int_{\Om\setminus K} |\nabla u|^p \, dxdy   & 
\mbox{ if } u \, \in \, W^{1,p}(\Om \setminus K) \mbox{ and }u=g \mbox{ on }\partial_D \Om,\\ 
\salt
\displaystyle +\infty & \mbox{ otherwise}.\end{cases}$$

Hence in this case the conclusion of Theorem \ref{mainr} follows from a general result
on convergence of minima (see Proposition \ref{Gamma-conv-prop}).

\begin{proof} {\bf  (i) $\Gamma$-liminf}: Let $u_h\to u$ strongly in  $L^p(\Om)$, 
we want to prove that $\liminf_{h\to\infty}F_h(u_h)\,\geq\,F(u)$. 
We can assume that $\liminf_{h\to\infty}F_h(u_h)=\lim_{h\to\infty}F_h(u_h)<\infty$. 
So,  for any $\Om'\subset\subset \Om\setminus K$ with $\partial_D\Om\subset\partial\Om'$, we have that 
$u_h\in W^{1,p}(\Om')$ for $h$ big enough, $u_h \rightharpoonup u$ in $W^{1,p}(\Om')$ and 
  $u=g$ on $\partial_D \Om $. 
Now from the lower semicontinuity of the
$L^p-$norm of the gradients and from the arbitrariness of $\Om'$ we get that $u\in W^{1,p}(\Om\setminus K)$
 and the $\Gamma$-liminf inequality holds.\\
{\bf (ii) $\Gamma$-limsup}:  
Let $u \,\in\, L^p(\Om)$. We want to construct a sequence $(u_h)\subset L^p(\Om)$ 
converging strongly to $u$ in $L^p(\Om)$  such that
$\lim_h F_h(u_h) \leq F(u)$. We can assume that $u\in W^{1,p}(\Om\setminus K)$ and $u=g$ on $\partial_D \Om $. 
 We set $u_h:=u$ in $\Om\setminus R_h$ 
where $R_h:= \left(-1,1\right)\times \left[ -\frac{1}{h},\frac{1}{h}\right]$.
Now let us define the function $u_h$  in $R_h$.  To this aim, we consider the function $v_h$  defined in $R_h$ 
by  $$v_h(x,y):= u\Bigl(x,\frac{2}{h} sgn(x) -y \Bigr)$$
where $sgn(x)$ denotes the sign of $x$. 
In other words, the function $v_h$ is obtained from $u$ by symmetry with respect to the segment $[0,1]\times\{\frac{1}{h}\}$ for $x$
positive and by symmetry with respect to the segment  $[-1,0]\times\{-\frac{1}{h}\}$ for $x$ negative.

Such a function $v_h$ may jump on the segment $\{0\} \times \left[-\frac{1}{h},\frac{1}{h}\right]$.
So, we consider the function  
 $\var \in C^0(R_h)$  defined by
$$
\var (x,y):=\begin{cases}
\,\,\,1 & \mbox{ if } |x| > \frac{1}{2},\\
-2x & \mbox{ if } - \frac{1}{2}\leq x \leq 0,\\
\,\,\,2x     & \mbox{ if }             0\leq  x \leq \frac{1}{2}.
\end{cases}
$$
Now we set $u_h:= \var \,v_h$ on $R_h$.
For this choice of $u_h$, it is easy to see that $u_h\in W^{1,p}(\Om\setminus K_h)$ with $u_h=g$ on $\partial_D \Om$ and that
the   $\Gamma$-limsup inequality holds.
\end{proof}
}
\end{example}

\begin{center}
\psfig{figure=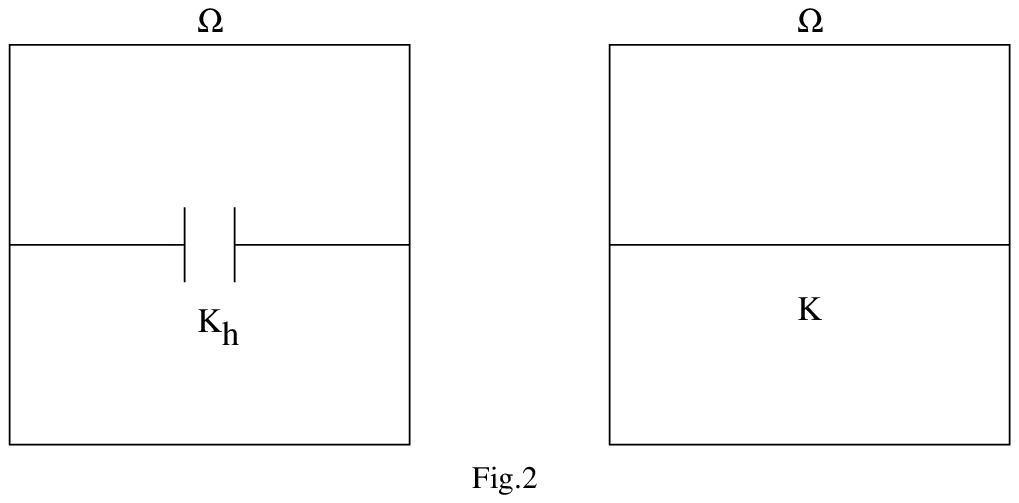}
\end{center}

In the following example, we consider a sequence of compact sets $K_h$ along which the problem $(P)$ is
not stable. More precisely in the limit problem, that is the problem solved by the limit function $u$,
there is an additional term involving the jump of $u$ on a point of $K$.  

\begin{example}\label{ex1}
{\rm Let $\Om$, $\partial _D\Om$ and $K$  be as in the previous example and let
$$K_h:=\left [-1,-\frac{a_h}{2}\right ]\times\Bigl\{0\Bigr\}
\cup\left [\frac{a_h}{2},1\right]\times\Bigl\{0\Bigr\}
\cup\left\{-\frac{a_h}{2}\right\}\times\left [-\frac{b_h}{2},\frac{b_h}{2}\right]
\cup\left\{\frac{a_h}{2}\right\}\times\left [-\frac{b_h}{2},\frac{b_h}{2}\right]$$
be as in  Fig. $2$ with $(a_h)$ and $(b_h)$ being two sequences of positive numbers converging to $0$. In this way $(K_h)$ converges to $K$ in the Hausdorff metric.
Let $F_h$ be  defined as in (\ref{gammasuc}). 

Assume that the sequence  $(\frac{1}{p} \,a_h \, b_h^{1-p})$ converges to some $ c\,\in\, [0,+\infty]$.
Then $F_h$ $\Gamma$-converges in the strong topology of $L^p(\Om)$ to 
$F_\infty$ defined in $L^p(\Om)$ in the following way (with the convention that $0  \cdot  \infty \,=\,0$).
\begin{equation}\label{gammalim}
F_\infty (u):= \begin{cases}
\displaystyle\frac{1}{p} \int_{\Om\setminus K} |\nabla u|^p dxdy  + 
c\, \left|u^+(0,0) - u^-(0,0)\right|^p\ & \mbox{ if }\left\{ \begin{array}{ll}
u \in  W^{1,p}(\Om \setminus K) \mbox{ and } &\mbox{ }\\
 u=g \mbox{ on }\partial_D \Om, & \mbox{ }
\end{array}
\right.\\ 
\salt
\displaystyle +\infty & \mbox{ otherwise, } \end{cases}\end{equation}
where $u^+(0,0)$ and $u^-(0,0)$ are respectively the values in $(0,0)$ of the traces of $u|_{\Om^+}$ and  $u|_{\Om^-}$ on $K$, 
$\Om^+$ and $\Om^-$ being respectively the upper and the lower connected components of $\Om \setminus K$.
\begin{proof}
 {\bf  (i) $\Gamma$-liminf}: Let $u_h\to u$ in  $L^p(\Om)$, we want to prove that 
$\liminf_{h\to\infty}F_h(u_h)\,\geq\,F(u)$. 
We can assume that $\liminf_{h\to\infty}F_h(u_h)=\lim_{h\to\infty}F_h(u_h)<\infty$. 
So,  for any $\Om'\subset\subset \Om\setminus K$ with $\partial_D\Om\subset\partial\Om'$, we have that 
$u_h\in W^{1,p}(\Om')$ for $h$ big enough, $u_h \rightharpoonup u$ in $W^{1,p}(\Om')$ and 
  $u=g$ on $\partial_D \Om $. 
Now from the lower semicontinuity of the
$L^p-$norm of the gradients and from the arbitrariness of $\Om'$ we get that $u\in W^{1,p}(\Om\setminus K).$

We set $R_h:=\Bigl(-\frac{a_h}{2}\,,\,\frac{a_h}{2}\Bigr)\times\Bigl(-\frac{b_h}{2}\,,\,\frac{b_h}{2}\Bigr).$ We have
\begin{eqnarray}\label{fe0}
\nonumber F_h(u_h)&=&\frac{1}{p} \int_{\Om\setminus K_h} |\nabla u_h|^p \, dxdy
\,=\, \frac{1}{p} \int_{\Om\setminus R_h} |\nabla u_h|^p \, dxdy\,+\,\frac{1}{p} \int_{R_h} |\nabla u_h|^p \, dxdy\\
\salt
& \geq &\frac{1}{p} \int_{\Om\setminus R_h} |\nabla u_h|^p \, dxdy \,+
\,\frac{1}{p}a_hb_h^{1-p}\intcauchy _{-\frac{a_h}{2}}^{\frac{a_h}{2}}\left|u_h\Bigl(x,-\frac{b_h}{2}\Bigr)\,-\,u_h\Bigl(x,\frac{b_h}{2}\Bigr)\right|^p dx.
\end{eqnarray}
 
Now let us fix  $\Om '\subset\subset \Om \setminus K$. We have that 
$\Om '\subset\subset \Om \setminus R_h$ definitively, and 
$$\liminf_{h\to\infty} \int_{\Om\setminus R_h} |\nabla u_h|^p \, dxdy \,\geq\,
\liminf_{h\to\infty} \int_{\Om'} |\nabla u_h|^p \, dxdy \,\geq\,\int_{\Om'} |\nabla u|^p \, dxdy.$$
By the arbitrariness of $\Om'$, we get

\begin{equation}\label{fe1}
\liminf_{h\to\infty} \int_{\Om\setminus R_h} |\nabla u_h|^p \, dxdy \,\geq\,\int_{\Om\setminus K} |\nabla u|^p \, dxdy.
\end{equation}
Let us consider  the functions $\tilde{u}_h^1$ defined in $(-1,1)\times(0,1- \frac{b_h}{2})$ by
$$\tilde{u}_h^1(x,y):= u_h|_{(-1,1)\times (\frac{b_h}{2},1)}\bigl(x,y + \frac{b_h}{2}\bigr)$$
and   $\tilde{u}_h^2$ defined in $(-1,1)\times(-1 + \frac{b_h}{2},0)$ by
 $$\tilde{u}_h^2(x,y):= u_h|_{(-1,1)\times (-1, - \frac{b_h}{2})}\bigl(x,y- \frac{b_h}{2}\bigr).$$
We  extend $\tilde{u}_h^1$ and $\tilde{u}_h^2$ respectively in $\Om^+$ and $\Om^-$ 
in such a way those extensions converge weakly
to $u$ respectively in $W^{1,p}(\Om^+)$ and  in $W^{1,p}(\Om^-)$. Recalling that $p>2,$
We have the uniform convergence of their traces on $K$.
So, 
\begin{eqnarray*}
\frac{1}{p}a_hb_h^{1-p}\intcauchy _{-\frac{a_h}{2}}^{\frac{a_h}{2}}\left|u_h\Bigl(x,-\frac{b_h}{2}\Bigr)
\,-\,u_h\Bigl(x,\frac{b_h}{2}\Bigr)\right|^p\,dx =
\frac{1}{p}a_hb_h^{1-p}\intcauchy _{-\frac{a_h}{2}}^{\frac{a_h}{2}}\left|\tilde u^1_h(x,0)
\,-\,\tilde u^2_h (x,0)\right|^p\,dx \\
\salt
 = \frac{1}{p}a_hb_h^{1-p}\intcauchy _{-\frac{a_h}{2}}^{\frac{a_h}{2}}\left|u^+(x,0)\,-\,u^-(x,0) + w_h(x) \right|^p\,dx,\end{eqnarray*}
with $(w_h)$ converging uniformly to $0$ on $K$. From this, it follows that
\begin{equation}\label{fe2} 
\lim_{h\to\infty} \frac{1}{p}a_hb_h^{1-p}\intcauchy _{-\frac{a_h}{2}}^{\frac{a_h}{2}}\left|u_h\Bigl(x,-\frac{b_h}{2}\Bigr)
\,-\,u_h\Bigl(x,\frac{b_h}{2}\Bigr)\right|^p\,dx \,=\, c \left|u^+(0,0)-u^-(0,0)\right|^p.
\end{equation}
Therefore, the  $\Gamma$-liminf inequality follows from (\ref{fe0}), (\ref{fe1}) and  (\ref{fe2}).\\
{\bf (ii) $\Gamma$-limsup}:  
Let $u \,\in\, L^p(\Om)$. We want to construct a sequence $(u_h)\subset L^p(\Om)$ which converges to $u$ such that
$\lim_h F_h(u_h) \leq F(u)$. We can assume that $u\in W^{1,p}(\Om\setminus K)$ and $u=g$ on $\partial_D \Om $. 

We set  $u_h=u$ in $(\Om\setminus K)\setminus R_h$ and we modify suitably $u$ in $R_h\setminus K$ 
in order to get a new function which does not jump on $K\cap R_h$. To this aim let $R_h^1:=R_h\cap\{y>0\}$ and $R_h^2:=R_h\cap\{y<0\}$. 
 Let us define $u_h$  in $R_h^1$. We set
$$v_h:=u|_{\bigl(-\frac{a_h}{2},\frac{a_h}{2}\bigr)\times\bigl(\frac{b_h}{2},b_h\bigr)}$$
and $$\tilde u_h(x,y):=v_h(x,b_h-y)\quad\mbox{ for any }(x,y)\in R^1_h.$$
In other words, $\tilde u_h$ is defined by taking the reflection of the restriction of $u$ on the 
rectangle symmetric to $R^1_h$ with respect to the horizontal line $y=\frac{b_h}{2}$. 
Now we consider the linear function  $\var_1(x,y):=\frac{2}{b_h}\,y.$ For any $(x,y)\in R^1_h$ we set
$$u_h(x,y):=\var_1(x,y)\Bigl(\tilde u_h(x,y)\,-\,\frac{u^+(0,0)+u^-(0,0)}{2}\Bigr)\,+\,\frac{u^+(0,0)+u^-(0,0)}{2}.$$
In the similar way, we define $u_h$ in $R_h^2$ using
$$u|_{\bigl(-\frac{a_h}{2},\frac{a_h}{2}\bigr)\times\bigl(-b_h,-\frac{b_h}{2}\bigr)}
\quad\mbox{ and }\quad\var_2(x,y):= -\frac{2}{b_h}\,y.$$
It is easy to check that $u_h\in W^{1,p}(\Om\setminus K_h)$, 
$u_h=g$ on $\partial_D \Om $ and by construction
$$\lim_{h\to\infty}\frac{1}{p}\int_{R_h^1} |\nabla u_h|^p \, dxdy\,=\,
\lim_{h\to\infty}\frac{1}{p}\int_{R_h^1} |\nabla u_h|^p \, dxdy\,=\,\frac{c}{2}\left|u^+(0,0)-u^-(0,0)\right|^p.$$
Therefore,
\begin{eqnarray*}
\lim_h F_h(u_h)&=&
\lim_h \frac{1}{p} \int_{\Om\setminus K_h} |\nabla u_h|^p \, dxdy\\
\salt
 &=&\lim_h  \frac{1}{p} \int_{(\Om\setminus K)\setminus R_h} |\nabla u_h|^p \,
dxdy\,+ \lim_h
 \,\frac{1}{p} \int_{R_h^1} |\nabla u_h|^p \, dxdy\,+\, \lim_h
\frac{1}{p} \int_{R_h^2} |\nabla u_h|^p \, dxdy\\ &=&
 \frac{1}{p}
\int_{\Om\setminus K} |\nabla u|^p \, dxdy\,
+\,c\left|u^+(0,0)-u^-(0,0)\right|^p\,=\,F_\infty(u).
 \end{eqnarray*}

\end{proof}
}
\end{example}
\begin{remark}
{\rm
Note that if the constant $c$ in the previous example is equal to zero, then 
we have the stability of $K$ for the minimization problem (\ref{eq1}) along the
sequence $K_h$ even if the intersection
 of the limit of the two connected
components of $K_h$ is the point 
 $(0,0)$ whose $(1,q)-$capacity is equal to
zero (recall that $q<2$). So, in Theorem \ref{mainr}
 the assumption that the
limit of   two connected components of $K_h$ is either empty or has positive
$(1,q)-$capacity is not a necessary condition. Although, it can not be removed
as shown by the case $c>0$.
 }
 \end{remark}
 \begin{remark}\label{tretracce}
{\rm
Starting from Example \ref{ex1}, one can construct examples in which the $\Gamma$-limit involves  
traces at the origin from more than two subdomains, as shown   in fig. $3$.   

\begin{center}
\psfig{figure=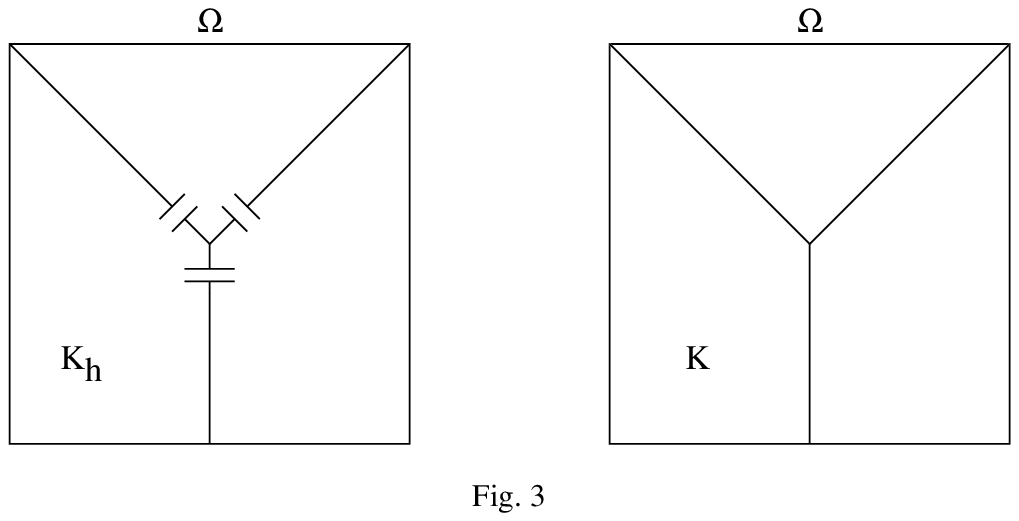}
\end{center}

In this case, we can obtain  a $\Gamma$-limit of the form 
$$F_\infty (u):= \begin{cases}
\displaystyle\frac{1}{p} \int_{\Om\setminus K} |\nabla u|^p dxdy  + 
\sum_{1\leq i<j\leq 3} c_{i,j}\, \left|u^i(0,0) - 
u^j(0,0)\right|^p\ & \mbox{ if }\left\{\begin{array}{ll}
 u \in  W^{1,p}(\Om \setminus K)&{ }\\
u=g \mbox{ on }\partial\Om & { }
\end{array}\right.\\ 
 
\salt
\displaystyle +\infty & \mbox{ otherwise },  \end{cases}$$
where $u^i(0,0)$ is the value at $(0,0)$ of the trace of $u|_{\Om_i}$,  $\Om_i$ being
the connected components of $\Om\setminus K$}.
\end{remark}

In Theorem \ref{mainr}, the assumption that $\Om$ is 
simply connected cannot be removed.
In fact we will consider in the next example a sequence of 
connected compact sets $K_h$ converging to $K$  and along which the stability of $K$
for the problem  (\ref{eq1}) does not hold.

\begin{example}\label{ex5}
{\rm
Let $\Om: = Q_2\setminus\overline Q_1$ with $Q_1$ and  $Q_2$  as in fig. 4. and let $K_h $ and $K$ be as in fig. 4.

\begin{center}
\psfig{figure=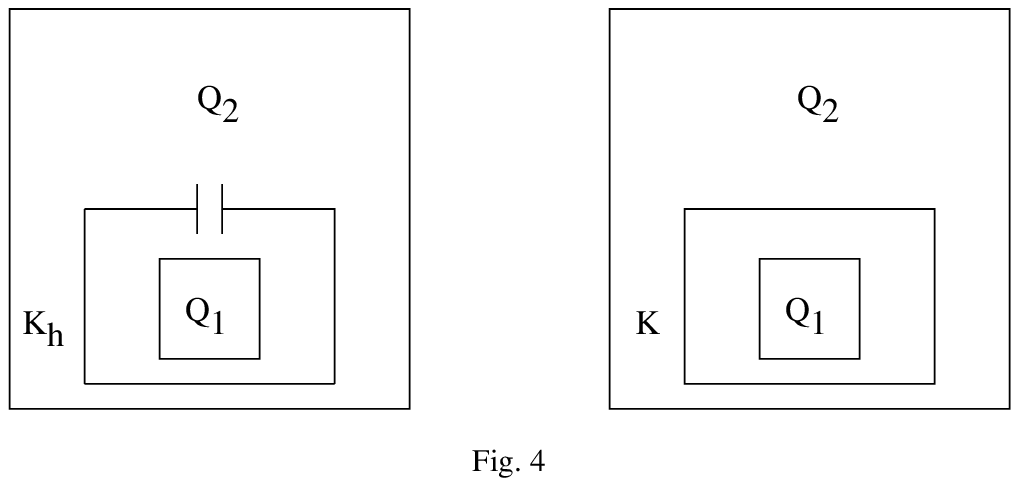}
\end{center}

In this case, arguing as in Example \ref{ex1} we have 
that $(F_h)\,\,\Gamma$-converge in the strong topology of $L^p(\Om)$ to 
the functional 
\begin{equation}\label{consc}
F_\infty (u):= \begin{cases}
\displaystyle\frac{1}{p} \int_{\Om\setminus K} |\nabla u|^p dxdy  + 
c\, \left|u^+(0,0) - u^-(0,0)\right|^p\ & \mbox{ if }\left\{ \begin{array}{ll}
u \in  W^{1,p}(\Om \setminus K) \mbox{ and } &\mbox{ }\\
 u=g \mbox{ on }\partial_D \Om, & \mbox{ }
\end{array}
\right.\\ 
\salt
\displaystyle +\infty & \mbox{ otherwise. } \end{cases}\end{equation}

Let $\partial_D\Om:= \partial\Om= \partial Q_1 \cup \partial Q_2$ and let 
$$g=\begin{cases}0 & \mbox{ on }  \partial Q_1,\\  1 & 
\mbox{ on }  \partial Q_2.\end{cases}$$ 
Let $u_h \in W^{1,p}(\Om\setminus K_h)$ be a solution of problem
 (\ref{eq1}) in  $\Om\setminus K_h$.
Then $\nabla u_h \rightharpoonup \nabla \tilde u$ 
weakly in $L^p(\Om,\R^2),$ where $\tilde u$ minimizes the functional $F_\infty$
among all functions $ w \in W^{1,p}(\Om\setminus K)$ with $w=g$ on $\partial_D\Om$. 
Now, let $\var \in C^1(\Om)$ with $\var=g$ on $\partial_D\Om$.
We can always assume that the sequence $(K_h)$ is such that  $c > \int_\Om |\nabla \var|^pdx$.
It is easy to see that the solution $u$ of Problem (\ref{eq1}) with data $g$
has gradient equal to $0$. 
 Then $$F_\infty (u) = c >  \int_\Om |\nabla
\var|^pdx = F_\infty(\var) \geq F_\infty (\tilde u).$$
 So, $u\neq \tilde u$
and hence $K$ is not stable along the sequence  $K_h$ 
 for Problem
(\ref{eq1}).  
  
}
\end{example}

\subsection{Limit problem obtained by duality}

In this section we examine, by a duality approach, 
the problem solved by the limit function $u$, even when there is not stability. 
Let $(K_h)\subset {\mathcal K}(\overline\Om)$ be such that $K_h$ converges 
to $K$ in the Hausdorff metric and  $| K_h|$ converges to $| K|$.
Let $u_h$ be  solution of (\ref{eq1}) in $(\Om\setminus K_h)$ and $v_h$
its conjugate which in this case satisfies
 \begin{equation}\label{pconj}
R\nabla v_h =|\nabla u_h|^{p-2}\nabla u_h\mbox{ a.e in }\Om.
\end{equation}
 From Lemma \ref{lemmaaux} it follows that, up to a subsequence, 
$\nabla u_h\rightharpoonup\nabla u$ weakly
 in $L^p(\Om,\R^2)$ for some function  $u\in L^{1,p}(\Om\setminus K)$ 
and $\nabla v_h\rightharpoonup\nabla v$ weakly
 in $L^q(\Om,\R^2)$ for some function  $v\in W^{1,q}(\Om)$. Using the fact that 
for every $\Om'\subset\subset\Om\setminus K$, 
we have ${\rm div}\,(|\nabla u_h|^{p-2}\nabla u_h)=0$ 
in ${\mathcal D}'(\Om')$ for $h$ big enough,
 it follows from the result in \cite{boccardo} that 
$\nabla u_h\to\nabla u$ a.e. in $\Om'$.  So by the arbitrariness of $\Om'$, 
 we get $\nabla u_h\to\nabla u$ a.e. in $\Om\setminus K$. Hence, 
using the fact 
$|K_h|\to|K|$  we can pass to the limit in (\ref{pconj}) and obtain 
\begin{equation}\label{pconj2}
R\nabla v =|\nabla u|^{p-2}\nabla u\mbox{ a.e in }\Om,
\end{equation}
through which we will find the limit problem solved
 by the function $u$ in the next example.
\vskip .3truecm
To this aim, we  call  contact point, any point 
of $K\cup\partial_N\Om$ which is limit of at least 
two sequences belonging to two different connected
 components of $K_h\cup\partial_N\Om$. 

\begin{center}
\psfig{figure=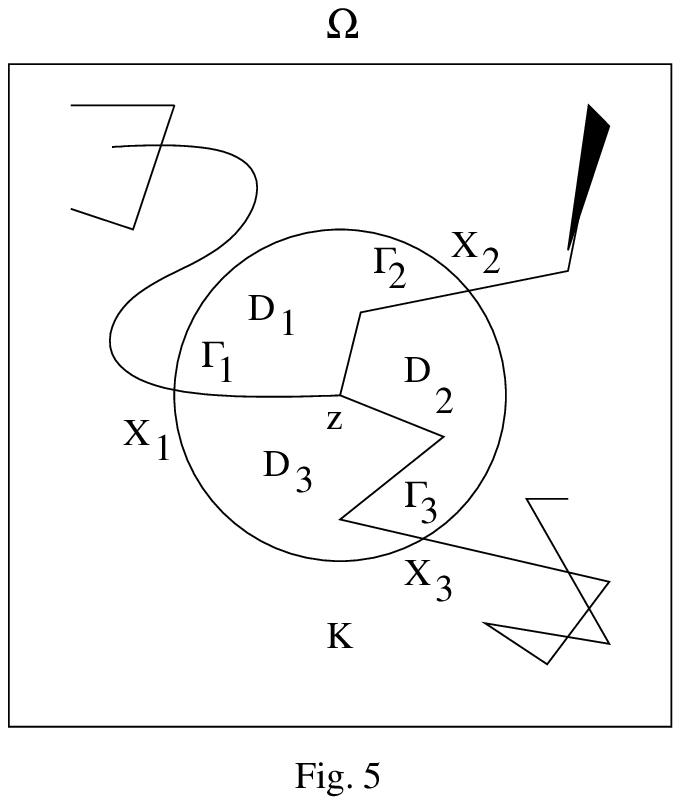}
\end{center}

\begin{example}\label{exfine}
 {\rm Let $(K_h)\subset {\mathcal K}_m(\overline\Om)$ be a sequence
 which converges to a compact set $K$ in the Hausdorff metric and such that 
 $| K_h|$ converges to $| K|$ with $K$ having only one contact point $z\in\Om$. 
Assume that there exists $r>0$ such that  
$B_{r}(z)\cap K=\Gamma_1\cup\Gamma_2\cup\Gamma_3$, 
with $\Gamma_i$ Lipschitz simple curves such
that $\Gamma_i\cap\Gamma_j = z $ for $i\neq j$ and 
$\Gamma_i \cap \partial B_r(z) = x_i$ for every $i$ 
(see  fig.$5$).
Suppose that $B_r(z)\setminus K=\bigcup_{i=1}^3 D_i$ with $D_i$ Lipschitz domains.  

If $u_h$ is a solution of (\ref{eq1}), then the sequence $(\nabla u_h)$
converges weakly to $\nabla u$ in $L^p(\Om,\R^2)$, 
 where the function $u$
solves a minimization problem of the type:
 \begin{equation}\label{minlim}
\displaystyle\min_{w}\left\{\begin{array}{ll}
\displaystyle\frac{1}{p}\int_{\Om\setminus K}|\nabla w|^pdx-
\bigl[a_1(w_3 -w_1)+a_2(w_1 -w_2)+a_3(w_2 -w_3)\bigr]\\
\salt
\hskip 4truecm w\in L^{1,p}(\Om\setminus K),\, 
w=g\mbox{ on }\partial_D\Om\setminus K & \hskip -.4truecm
\end{array}
\right\}
\end{equation}
where $w_j$ is the trace of $w|_{D_j}$ evaluated at $z$ and $a_j$ coincides with the value taken 
on $\Gamma_j$ by the continuous representative of 
 the limit $v$ of the conjugates $v_h$.}

\begin{proof}

First of all from Lemma \ref{lemmaaux} it follows that
 $u\in L^{1,p}(\Om\setminus K)$ and $u=g$ on $\partial_D\Om\setminus K$. 
Now let $\var\in C^{1}(\Om\setminus K)\cap L^{1,p}(\Om\setminus K)$ with 
$\var|_{D_j}\in C^1 (\overline{D_j})\,\, \forall j$ and  such that  $\var=0\,
\mbox{ on }\,\partial_D\Om\setminus K$. 
Using (\ref{pconj2}), we have that
\begin{eqnarray}\label{form0}
\nonumber \int_{\Om\setminus K}|\nabla u|^{p-2}\nabla u\nabla\var\,dx
 &=&\int_{\Om\setminus K}R \nabla v\,\nabla\var\,dx\\
&=&\int_{(\Om\setminus K)\setminus B_{r}(z)}R \nabla  v\,\nabla\var\,dx\,+
\, \sum_{j=1}^3\int_{D_j}R\nabla v \nabla\var\,dx
\end{eqnarray}
Since $z$ is the only contact point of $K\cup\partial_N\Om$, we have
 that any connected component of
\mbox{$K\cup\partial_N\Om \setminus B_r(z)$} is contained in a limit of some 
connected component of $ K_h\cup\partial_N\Om.$ Therefore
 $v$ is constant on every connected components of 
$K\cup\partial_N\Om \setminus B_r(z)$. In order to
integrate by parts outside $B_r(z)$, where $K$ in general is not regular,  
it is useful (from \cite[Theorem 4.5]{hekima}) to  approximate strongly in 
$W^{1,q}(\Om\setminus \overline{B_r(z)})$ the function $ v$
with smooth functions $w_n$ that are constant in a suitable 
neighborhood of any connected component of 
$(K\cup \partial_N\Om)\setminus B_r(z)$. As $ \text{div}\,(R\nabla w_n)=0$,
integrating by parts we get
\begin{eqnarray}\label{form1}
\nonumber  & &\int_{(\Om\setminus K)\setminus\overline{B_{r}(z)}}R \nabla v\,
\nabla\var\,dx \,\,= \,\,\lim_{n\to\infty} 
\int_{(\Om\setminus K)\setminus B_{r}(z)}R \nabla w_n\,\nabla\var\,dx\\
\nonumber & & \,\,=\,\,\- \lim_{n\to\infty}
\int_{\partial B_r(z)} (R \nabla w_n)\nu\,\var\,d\mathcal H^1\,\,=\,\,-
\lim_{n\to\infty} \int_{\partial B_r(z)}  
\frac{\partial w_n}{\partial\tau} \var  \, d\mathcal H^1 \\
\salt
\nonumber & &\,\, =\,\, + \int_{\partial B_r(z)} 
v \frac{\partial\var}{\partial\tau} \,d\mathcal H^1 -(a_1 \var|_{D_1}(x_1) 
- a_2 \var|_{D_1}(x_2))
- (a_2 \var|_{D_2}(x_2) - a_3 \var|_{D_2}(x_3)) \\
\salt
& & \qquad - (a_3 \var|_{D_3}(x_3) - a_1 \var|_{D_3}(x_1))
\end{eqnarray}
where $\nu$ is the unit vector outer normal to $B_r(z)$ and $\tau= -R \nu$ is 
the corresponding tangential unit vector, so that $\frac{\partial\var}{\partial\tau}$
 denotes the tangential derivative of the function $\var$.

On the other hand, as $\text{div}\,(R\nabla\var)=0$, we have
\begin{equation}\label{eqinter1}
\sum_{j=1}^3 \int_{D_j} R \nabla v\nabla\var\,dx \,=\,
- \sum_{j=1}^3 \int_{D_j}  \nabla v R \nabla\var\,dx \,=\,
- \sum_{j=1}^3 \int_{\partial D_j} v \frac{\partial\var}{\partial\tau} \,d\mathcal H^1.
\end{equation}
Let us compute 
 \begin{eqnarray*}
 & & \int_{\partial D_1} v \frac{\partial\var}{\partial\tau} \,d\mathcal H^1 \,\, =\,\,
\int_{\partial B_r(z)\cap\partial D_1} v \frac{\partial\var}{\partial\tau} \,d\mathcal H^1
+\int_{\Gamma_1} v \frac{\partial\var}{\partial\tau} \,d\mathcal H^1 +\int_{\Gamma_2}
 v \frac{\partial\var}{\partial\tau} \,d\mathcal H^1\\
\salt
& = &  \int_{\partial B_r(z)\cap\partial D_1} v \frac{\partial\var}{\partial\tau}
 \,d\mathcal H^1 - a_1(\var|_{D_1}(x_1) - \var|_{D_1}(z)) - a_2
(\var|_{D_1}(z) - \var|_{D_1}(x_2)).
\end{eqnarray*}

In a similar way
$$
 \int_{\partial D_2} v \frac{\partial\var}{\partial\tau} \,d\mathcal H^1 =
  \int_{\partial B_r(z)\cap\partial D_2} v \frac{\partial\var}{\partial\tau}
 \,d\mathcal H^1 - a_2(\var|_{D_2}(x_2) - \var|_{D_2}(z)) - a_3
(\var|_{D_2}(z) - \var|_{D_2}(x_3))$$

$$ \int_{\partial D_3} v \frac{\partial\var}{\partial\tau}
 \,d\mathcal H^1= 
 \int_{\partial B_r(z)\cap\partial D_3} v 
\frac{\partial\var}{\partial\tau} \,d\mathcal
 H^1 -a_3(\var|_{D_3}(x_3) - \var|_{D_3}(z)) - a_1
(\var|_{D_3}(z) - \var|_{D_3}(x_1)).
$$

Then, we have by (\ref{eqinter1}) that
\begin{eqnarray}\label{form2}
\nonumber & &
\sum_{j=1}^3 \int_{D_j} R \nabla v\nabla\var\,dx \,=\, 
 -\int_{\partial B_r(z)} v \frac{\partial\var}{\partial\tau} \,d\mathcal H^1
 + a_1(\var|_{D_1}(x_1) - \var|_{D_1}(z))\\
\salt
\nonumber & &  + a_2(\var|_{D_1}(z) - \var|_{D_1}(x_2)) +
 a_2(\var|_{D_2}(x_2) - \var|_{D_2}(z)) + a_3(\var|_{D_2}(z) - \var|_{D_2}(x_3)) \\
\salt
\nonumber & & + a_3(\var|_{D_3}(x_3) - \var|_{D_3}(z)) +
 a_1(\var|_{D_3}(z) - \var|_{D_3}(x_1))\\
\salt
\nonumber & & = -\int_{\partial B_r(z)} v \frac{\partial\var}{\partial\tau} 
\,d\mathcal H^1 +a_1(\var|_{D_3}(z)-\var|_{D_1}(z))
+a_2(\var|_{D_1}(z)-\var|_{D_2}(z))\\
\salt
\nonumber & & +a_3(\var|_{D_2}(z)-\var|_{D_3}(z)) +
\left[(a_1 \var|_{D_1}(x_1) - a_2 \var|_{D_1}(x_2))
+ (a_2 \var|_{D_2}(x_2) - a_3 \var|_{D_2}(x_3))\right. \\
\salt
& & \qquad \left.+ (a_3 \var|_{D_3}(x_3) - a_1 \var|_{D_3}(x_1))\right].
\end{eqnarray}
Therefore from (\ref{form0}), (\ref{form1}) and (\ref{form2}) we get the identity
$$
\int_{\Om\setminus K}|\nabla u|^{p-2}\nabla u\nabla\var\,
dx=a_1(\var|_{D_3}(z)-\var|_{D_1}(z))
+a_2(\var|_{D_1}(z)-\var|_{D_2}(z))+a_3(\var|_{D_2}(z)-\var|_{D_3}(z))
$$
which is the weak formulation of the Euler-Lagrange equation 
for the minimization problem (\ref{minlim}) 
and the conclusion follows from the convexity of the functional in  (\ref{minlim}).
\end{proof}
\end{example}
   
\begin{remark}\label{remai}
{\rm Note that the function $u$ solves a problem of the type (\ref{minlim}) 
whenever $K$ is as in Fig.5 and $z$ is the only contact point of $K$, 
independently of the sequence $(K_h)$. 
However the constants $a_i$ are related to the particular 
sequence $(K_h)$ and to the particular boundary data $g$. 
Indeed,  the constants $a_i$ are the limits (as $h\to\infty$) of the values 
taken by the conjugates $v_h$ on the connected components of $K_h$. 

Note also that the functional to minimize in (\ref{minlim}) can be rewritten as
$$\displaystyle\frac{1}{p}\int_{\Om\setminus K}|\nabla w|^pdx-
\bigl[w_1(a_2-a_1)+ w_2(a_3-a_2)+w_3(a_1-a_3)\bigr].$$
}
\end{remark}
\begin{remark}\label{remgamconj}
{\rm The example above can easily be extended to  cases in which 
there are finitely many 
contact points where a finite number of curves intersect each 
other as above. In particular, this method applied to Example \ref{ex1} 
gives the following
 minimization problem in the limit
\begin{equation}\label{minlim2}
\displaystyle\min_{w}\left\{\begin{array}{ll}
\displaystyle\frac{1}{p}\int_{\Om\setminus K}|\nabla w|^pdx-
(a_1-a_2)(u^+(0,0)-u^-(0,0))\\
\salt
\hskip 4truecm w\in L^{1,p}(\Om\setminus K),\, 
w=g\mbox{ on }\partial_D\Om\setminus K & \hskip -.4truecm
\end{array}
\right\}.
\end{equation}
Now from the weak Euler-Lagrange equations of (\ref{minlim2}) and 
of the minimization problem involving the functional $F_\infty$ in 
(\ref{gammalim})
we get  for every $\var\in L^{1,p}(\Om\setminus K)$ with $\var=0$ on 
$\partial_D\Om\setminus K$
\begin{eqnarray*}
& & pc|u^+(0,0)-u^-(0,0)|^{p-2}(u^+(0,0)-u^-(0,0))(\var^+(0,0)-\var^-(0,0))\,=\\
\salt
& &\hskip 3truecm  =\,\,(v^+(0,0)-v^-(0,0))(\var^+(0,0)-\var^-(0,0)).
\end{eqnarray*}
 By the arbitrariness of $\var$ we get 
$$pc|u^+(0,0)-u^-(0,0)|^{p-2}(u^+(0,0)-u^-(0,0))=(v^+(0,0)-v^-(0,0))$$
which can be interpreted as a discrete version of the duality relation (\ref{pconj2}).

}
\end{remark}
\vskip .5truecm
\centerline{\sc Acknowledgements}
 \vskip .1truecm
\noindent The authors wish to thank Gianni Dal Maso for having proposed 
the subject of this paper, and for many interesting discussions. 
This work is part of the
European Research Training Network ``Homogenization and Multiple
Scales'' under contract HPRN-2000-00109, and of the Research \
Project ``Calculus of Variations''
supported by SISSA and by
the Italian Ministry of Education, University, and Research.
 
\end{document}